\documentclass[12pt,twoside]{article}
\usepackage[ansinew]{inputenc} 
\usepackage{natbib}
\usepackage{latexsym,amsmath}
\usepackage{amsmath,amsthm}
\usepackage[psamsfonts]{amssymb}
\usepackage{amsmath}
\usepackage{amsfonts}

\usepackage{psfig, graphicx} 

\def\R{I\!\!R}

\def \be{\begin{equation}}
\def \ee{\end{equation}}

\makeatletter
\@addtoreset{equation}{section}

\makeatother
\newtheorem{theo}{Theorem}

\begin{document}

\begin{center}
{\Large \textbf{Uniform in bandwidth consistency of local polynomial
regression function estimators}}\\[0.5cm]{\large Julia Dony \\
Departement Wiskunde, Vrije Universiteit Brussel}\\{\large Uwe Einmahl
\\ Departement Wiskunde, Vrije Universiteit Brussel}\\%
{\large David M. Mason \\ Statistics Program, University of Delaware}
\end{center}

\begin{quote}
\textbf{Abstract:} {We generalize a method for proving uniform in bandwidth
consistency results for kernel type estimators developed by the two last named
authors. Such results are shown to be useful in establishing consistency of
local polynomial estimators of the regression function. }

\end{quote}

\section{Introduction}

\nocite{2}
\nocite{3}
Let $X,X_{1},X_{2},...$ be i.i.d. $\mathrm{I\!R}^{d}$  ($d\geq1)$ valued
random variables and assume that the common distribution function of these
variables has a Lebesgue density function, which we shall denote by $f_{X}$. A
kernel $K$ will be any measurable function which satisfies the following
conditions:
$$
\int_{\mathrm{I\!R}^{d}}K(s)ds=1, \quad \textrm{and} \leqno(K.i)
$$
$$
\Vert K\Vert_{\infty}:=\sup_{x\in\mathrm{I\!R}^{d}}|K(x)|=\kappa
<\infty.\leqno(K.ii)
$$
The kernel density estimator of $f_{X}$ based upon the sample $X_{1}%
,...,X_{n}$ and bandwidth $0<h<1$ is
$$
\widehat{f}_{n,h}(x):=\frac 1{nh}\sum_{i=1}^{n}K\left(\frac{x-X_{i}}{h^{1/d}}\right),\mbox{ }x\in
\mathrm{I\!R}^{d}.
$$
It is well known that if one chooses a suitable bandwidth sequence
$h_{n}\rightarrow0$ and the density $f_{X}$ is continuous, one obtains a
strongly consistent estimator $\widehat{f}_{n}:=\widehat{f}_{n,h_{n}}$ of
$f_{X}$, i.e. one has with probability $1$, $\widehat{f}_{n}(x)\rightarrow
f_{X}(x),x\in\mathrm{I\!R}^{d}$. It is also natural to investigate other modes
of convergence, for instance uniform convergence and to ask what convergence
rates are feasible.\newline

For proving such results, one usually writes the difference $\widehat{f}%
_{n}(x)-f_{X}(x)$ as the sum of a probabilistic term $\widehat{f}%
_{n}(x)-\mathrm{I\!E}\widehat{f}_{n}(x)$ and a deterministic term $\mathrm{I\!E}\widehat
{f}_{n}(x)-f_{X}(x),$ the so-called bias. The order of the bias depends on
smoothness properties of $f_{X}$ only, whereas the first (random) term can be
studied via empirical process techniques as has been pointed out by \cite{10}, \cite{11}, \cite{12} and \cite{9}, among other authors.\smallskip

\cite{7} (see also \cite{1b} for the 1-dimensional
case) have shown that if $K$ is a \textquotedblleft regular\textquotedblright%
\ kernel, the density function $f_{X}$ is bounded and $h_{n}$ satisfies the
regularity conditions $h_{n}\searrow0,$ $h_{n}/h_{2n}$ is bounded, and
$$
\log(1/h_{n})/\log\log n\rightarrow\infty\mbox{ and }nh_{n}/\log
n\rightarrow\infty,
$$
one has with probability $1$,
\begin{equation}
\Vert\widehat{f}_{n}-\mathrm{I\!E}\widehat{f}_{n}\Vert_{\infty}=O\left(  \sqrt{|\log
h_{n}|/nh_{n}}\right)  , \label{bb}%
\end{equation}
where $||\cdot||_{\infty}$ denotes the supremum norm on $\mathrm{I\!R}^{d}.$
Moreover, this rate cannot be improved. Interestingly one does not need
continuity of $f_{X}$ for this result. (Continuity of $f_{X}$ is of course
needed for controlling the bias.) \smallskip Recently, \cite{4a} have provided a \textquotedblleft uniform in $h$\textquotedblright%
\ version of this result, that is, they have proved that
\begin{equation}
\limsup_{n\rightarrow\infty}\sup_{\frac{c\log n}{n}\leq h\leq1}\frac{\sqrt
{nh}\Vert\widehat{f}_{n,h}-\mathrm{I\!E}\widehat{f}_{n,h}\Vert_{\infty}}{\sqrt{|\log
h_{n}| \vee\log\log n}}=:K(c)<\infty. \label{kd}%
\end{equation}
This result implies that if one chooses the bandwidth depending on the data
and/or the location $x$, as is usually done in practice, one has the same
order of convergence as in the case of a deterministic bandwidth sequence.

Now let $Y,Y_{1},Y_{2},\ldots$ be a sequence of $r$-dimensional random vectors
($r\geq1$) so that the random vectors $(X,Y),(X_{1},Y_{1}),\ldots$ are i.i.d.
with common joint Lebesgue density function $f.$ In this case it is also of
great interest to estimate $\mathrm{I\!E}\left[  \psi(Y)|X=x\right]  $, where
$\psi:I\!\!R^{r}\rightarrow I\!\!R$ is a suitable mapping. A possible kernel
type estimator which reduces to the classical Nadaraya-Watson estimator if
$r=1,\psi(y)=y$, is given by
\begin{equation}
\widehat{m}_{n}(x,\psi)=\frac{\sum_{i=1}^{n}\psi(Y_{i})K((x-X_{i})/h_{n}%
)}{\sum_{i=1}^{n}K((x-X_{i})/h_{n})}. \label{NW}%
\end{equation}
Likewise by setting in the 1-dimensional case for $t\in I\!\!R$, $\psi
_{t}(y)=I_{]-\infty,t]}(y),y\in I\!\!R$, we obtain the kernel estimator of the
conditional empirical function
$$
F(t|x):=I\!\!P\{Y\leq t|X=x\}
$$
given by
$$
\widehat{F}{_{n}}(t|x):=\frac{\sum_{i=1}^{n}1(Y_{i}\leq t)K((x-X_{i})/h_{n}%
)}{\sum_{i=1}^{n}K((x-X_{i})/h_{n})}.
$$
This kernel estimator is called the conditional empirical distribution
function and was first extensively studied by \cite{13}. Exact convergence
rates uniformly on compact subsets of $I\!\!R^{d}$ have been obtained for both
Nadaraya-Watson type estimators as in $(\ref{NW})$ and the conditional
empirical distribution function by \cite{4} in the case of
deterministic bandwidth sequences. Recently, \cite{4a} have
established uniform in bandwidth results for these estimators which are of a
similar type as result $(\ref{kd})$. The proof of these results requires
establishing a suitable version of a result of type $(\ref{kd})$ for processes
of the form
$$
\frac{1}{nh}\sum_{i=1}^{n}\left\{  \varphi(Y_{i})K\left(\frac{x-X_{i}}{h^{1/d}}
\right)-\mathrm{I\!E}\left[  \varphi(Y)K\left(\frac{x-X_{i}}{h^{1/d}}\right)\right]  \right\}  ,
$$
where $x\in I$ ($I$ a compact subset of $I\!\!R^{d}$ or $I=I\!\!R^{d}$) and
$\varphi\in\Phi$, where $\Phi$ is a suitable class of functions.\newline

For certain applications, however, this class of processes could be too small.
One of the purposes of this paper is to establish such uniform in bandwidth
consistency results for a larger class of processes. As an application of our
results, we shall prove uniform in bandwidth consistency of local polynomial
regression estimators. Such estimators are generalizations of the classic
Nadaraya-Watson estimator (see, especially, \cite{5} and
\cite{16}). In Section 2 we will state two general consistency results,
one of which will be proved in Section 3. In Section 4 we treat the local
polynomial regression estimators. In an appendix we gather together some facts
needed in our proofs.

\section{General consistency results}

We shall begin by stating a result proved in \cite{4a}, which
will be instrumental in establishing uniform in bandwidth consistency of local
polynomial regression function estimators. Let $\Phi$ denote a class of
measurable functions on $\mathrm{I\!R}^{r}$ with a finite valued measurable
envelope function $F$,
\begin{equation}
F(y)\geq\sup_{\varphi\in\Phi}|\varphi(y)|,\mbox{ }y\in\mathrm{I\!R}^{r}.
\label{ev}%
\end{equation}
Further assume that $\Phi$ is pointwise measurable and satisfies (A.2) in the
Appendix with $\mathcal{G}$ replaced by $\Phi.$ (For the definition of
pointwise measurable also refer to the Appendix.) Consider the following class
of functions
\begin{equation}
\mathcal{K}=\left\{  K((x-\cdot)/h^{1/d}):h>0,x\in\mathrm{I\!R}^{d}\right\}  ,
\label{KK}%
\end{equation}
and assume that $\mathcal{K}$ is pointwise measurable and satisfies (A.2) with
$\mathcal{G}$ replaced by $\mathcal{K}$. Introduce the class of continuous
functions on a compact subset $J$ of $\mathrm{I\!R}^{d}$ indexed by $\Phi$:%
$$
\mathcal{C}:=\{c_{\varphi}:\varphi\in\Phi\}.
$$
We shall always assume that the class $\mathcal{C}$ is relatively compact with
respect to the sup-norm topology, which by the Arzela-Ascoli theorem is
equivalent to being uniformly bounded and uniformly equicontinuous.\\

For any $\varphi\in\Phi$ and continuous functions $c_{\varphi}$ on a compact
subset $J$ of $\mathrm{I\!R}^{d}$, set for $x\in J,$
$$
\eta_{\varphi,n,h}(x)=\sum_{i=1}^{n}c_{\varphi}(x)\varphi(Y_{i})K\left(
\frac{x-X_{i}}{h^{1/d}}\right)  ,
$$
where $K$ is a kernel with \textit{support contained in} $\left[
-1/2,1/2\right]  ^{d}$ such that (K.i) and (K.ii) hold. The following result
was proved in \cite{4a}, where it is stated as Proposition 2.
($\left\Vert \cdot\right\Vert _{I}$ denotes the supremum norm on $I$.)\newline

\begin{theo}
Let $I$ be a compact subset of $\mathrm{I\!R}^{d}$ such that $J=I^{\eta},$ for
some $0<\eta<1$. Also assume that
\begin{equation}
f\mbox{ {\it   is continuous and strictly positive on} }J. \label{f3}%
\end{equation}
Further assume that the envelope function $F$ of the class $\Phi$ satisfies
\begin{equation}
\exists M>0:F(Y)1\left\{  X\in J\right\}  \leq M,\mbox{ a.s.} \label{mm}%
\end{equation}
or for some $p>2$
\begin{equation}
\alpha:=\sup_{z\in J}\mathrm{I\!E}[F^{p}(Y)|X=z]<\infty. \label{PP}%
\end{equation}
Then we have for any $c>0$ and $0<h_{0}<(2\eta)^{d},$ with probability 1,
\begin{equation}
\limsup_{n\rightarrow\infty}\sup_{c(\log n/n)^{\gamma}\leq h\leq h_{0}}%
\frac{\sup_{\varphi\in\Phi}\left\Vert \eta_{\varphi,n,h}-\mathrm{I\!E}\eta
_{\varphi,n,h}\right\Vert _{I}}{\sqrt{nh\left(  |\log h| \vee\log\log
n\right)  }}=:Q(c)<\infty, \label{PPP}%
\end{equation}
where $\gamma=1$ in the bounded case $(\ref{mm})$
and $\gamma=1-2/p$ under assumption $(\ref{PP})$.
\end{theo}

The next result generalizes Theorem 1 in the bounded case. Its proof is
illustrative of how that of Theorem 1 goes using an empirical process approach
based upon an inequality of Talagrand coupled with a moment bound for the
supremum of the empirical process. These basic tools are stated in the Appendix.\\

In the following, $||\cdot||_{\infty}$ denotes the supremum norm on
$\mathrm{I\!R}^{d}$ or $\mathrm{I\!R}^{d+r}$, whichever is appropriate. Let
$\mathcal{G}$ denote a class of measurable real valued functions $g$ of
$\left(  u,t\right)  \in\mathrm{I\!R}^{d}\times\mathrm{I\!R}^{r}%
=\mathrm{I\!R}^{d+r}$. We shall assume that $\mathcal{G}$
satisfies:\textit{\smallskip}

(G.i) $\sup_{g\in\mathcal{G}}$ $||g||_{\infty}=:\kappa<\infty;$%
\textit{\medskip}

(G.ii) $\sup_{g\in\mathcal{G}}\int_{\mathrm{I\!R}^{d+r}}g^{2}%
(x,y)dxdy=:L<\infty.$\textit{\medskip}

\noindent Denote by $\mathcal{F}_{\mathcal{G}}$, the class of functions of
$\left(  s,t\right)  $ $\in$ $\mathrm{I\!R}^{d+r}$ formed from $\mathcal{G}$
as follows\textit{:}
$$
\mathcal{F}_{\mathcal{G}}=\left\{  g(z-s\lambda,t):\lambda\geq1,z\in
\mathrm{I\!R}^{d}\text{ and }g\in\mathcal{G}\right\}  .
$$

\noindent We\ shall also assume that the class of functions $\mathcal{F}%
_{\mathcal{G}}$ satisfies the following uniform entropy
condition:\textit{\smallskip}

(F.i) for some $C_{0}>0$ and $\nu_{0}>0,$ $N(\epsilon,\mathcal{F}%
_{\mathcal{G}})\leq C_{0}\epsilon^{-\nu_{0}},$ $0<\epsilon<1.$%
\textit{\smallskip}

\noindent Finally, to avoid using outer probability measures in all of our
statements, we impose the measurability assumption:\textit{\smallskip}

(F.ii) $\mathcal{F}_{\mathcal{G}}$ is a pointwise measurable
class.\textit{\smallskip}

\noindent(For the definitions of pointwise measurable and of $N(\epsilon
,\mathcal{F}_{\mathcal{G}})$ see the Appendix below, where we use $\kappa$ as
our envelope function.)\textit{\smallskip}

For any $g\in\mathcal{G}$ and $0<h<1$ define,
$$
g_{n,h}(x):=(nh)^{-1}\sum_{i=1}^{n}g\left(\frac{x-X_{i}}{h^{1/d}},Y_{i}\right),\mbox{ }x\in
\mathrm{I\!R}^{d}.
$$
\textit{\smallskip}

\begin{theo} Assuming (G.i), (G.ii), (F.i), (F.ii),
and $f$ (the joint density of $(X,Y)$)  bounded, we have for
$c>0$ and $0<h_{0}<1$,
\begin{equation}
\limsup_{n\rightarrow\infty}\sup_{\frac{c\log n}{n}\leq h\leq h_{0}}\sup
_{g\in\mathcal{G}}\frac{\sqrt{nh}\left\Vert g_{n,h}-\mathrm{I\!E}%
g_{n,h}\right\Vert _{\infty}}{\sqrt{ |\log h| \vee\log\log n }}=:G(c)<\infty.
\label{gd}%
\end{equation}
\end{theo}

\paragraph{Remark.}

Theorem 2 is still valid for $r=0$. In this case, $g:\R^d\to\R^d$ and condition (G.ii) should be read as $\sup_{g\in\mathcal{G}}\int_{\mathrm{I\!R}^{d}}g^{2}%
(x)dx=:L<\infty$. 

\section{Proof of Theorem 2}

Let $\alpha_{n}$ be the empirical process based on the sample $\left(
X_{1},Y_{1}\right)  ,\ldots,\left(  X_{n},Y_{n}\right)  $, i.e. if
$\varphi:\mathrm{I\!R}^{d}\times\mathrm{I\!R}^{r}\rightarrow\mathrm{I\!R}$, we
have
$$
\alpha_{n}(\varphi)=\sum_{i=1}^{n}(\varphi(X_{i},Y_{i})-\mathrm{I\!E}\varphi
(X,Y))/\sqrt{n}.
$$
Notice that in this notation
$$
g_{n,h}\left(  x\right)  -\mathrm{I\!E}g_{n,h}\left(  x\right)  =\frac
{1}{h\sqrt{n}}\alpha_{n}\left(  g((x-\cdot)/h^{1/d},\cdot)\right), \quad x\in\R^d
$$
so we get that
$$
\sup_{g\in\mathcal{G}}\frac{\sqrt{nh}\left\Vert g_{n,h}-\mathrm{I\!E}%
g_{n,h}\right\Vert _{\infty}}{\sqrt{|\log h| \vee
\log\log n }}=\sup_{g\in\mathcal{G}}\sup_{x\in\mathrm{I\!R}^{d}}%
\frac{\left\vert \sqrt{n}\alpha_{n}\left(  g\left(  \frac{x-\cdot}{h^{1/d}%
},\cdot\right)  \right)  \right\vert }{\sqrt{nh\left(  |\log h|
\vee\log\log n\right)  }},
$$
where $g((x-\cdot)/h^{1/d},\cdot)$ denotes the function $\left(  s,t\right)
\rightarrow g((x-s)/h^{1/d},t)$. We first note that by (G.ii) and the
assumption that $||f||_{\infty}<\infty$,
\begin{align*}
\mathrm{I\!E}\left[  g^{2}\left(  \frac{x-X}{h^{1/d}},Y\right)  \right]   &
= h\int_{\mathrm{I\!R}^{d}}\int_{\mathrm{I\!R}^{r}}h^{-1}g^{2}\left(
\frac{x-s}{h^{1/d}},t\right)  f\left(  s,t\right)  dsdt\\
&  \leq h||f||_{\infty}L.
\end{align*}
Set for $j\geq0$ and $c>0,$%
$$
h_{j,n}:=\left(  2^{j}c\log n\right)  /n
$$
and
$$
{\mathcal{F}}_{j,n}=\left\{  g((x-\cdot)/h^{1/d},\cdot):g\in\mathcal{G}%
,h_{j,n}\leq h\leq h_{j+1,n},\text{ }x\in\mathrm{I\!R}^{d}\right\}  .
$$
Clearly for $h_{j,n}\leq h\leq h_{j+1,n},$%
$$
\mathrm{I\!E}\left[  g^{2}\left(  \frac{x-X}{h^{1/d}},Y\right)  \right]
\leq2h_{j,n}||f||_{\infty}L=:D_{0}h_{j,n}=:\sigma_{j,n}^{2}.
$$

We shall use Proposition A.1 in the Appendix to bound $\mathrm{I\!E}\Vert
\sum_{i=1}^{n}\epsilon_{i}\varphi(X_{i},Y_{i})\Vert_{{\mathcal{F}}_{j,n}}%
.$\textit{\ }To that end we note that each $\mathcal{F}_{j,n}$ satisfies (A.1)
of the proposition with $G=\beta=\kappa$ and (A.3) with $\sigma^{2}%
=\sigma_{j,n}^{2}$. Further, since $\mathcal{F}_{j,n}\subset\mathcal{F}%
_{\mathcal{G}},$ we see by (F.i) that each $\mathcal{F}_{j,n}$ also fulfills
(A.2). Finally (A.4) holds for large enough $n$ and all $j\geq0$. Now by
applying Proposition A.1 we get for all large enough $n$ and $j\geq0$,
\begin{equation}
\mathrm{I\!E}\Vert\sum_{i=1}^{n}\epsilon_{i}\varphi(X_{i},Y_{i})\Vert
_{\mathcal{F}_{j,n}}\leq D_{1}\sqrt{nh_{j,n}\left\vert \log\left(
D_{2}h_{j,n}\right)  \right\vert }, \label{II}%
\end{equation}
for some $D_{1}>0$ and $D_{2}>0$. Let for large enough $n$
$$
l_{n}:=\max\left\{  j:h_{j,n}\leq2h_{0}\right\}  ,
$$
then a little calculation shows that
\begin{equation}
l_{n}\sim\frac{\log\left(  \frac{nh_{0}}{c\log n}\right)  }{\log2}. \label{LL}%
\end{equation}
For $k\geq1$, set $n_{k}=2^{k}$, and let
$$
c_{j,k}:=\sqrt{n_{k}h_{j,n_{k}}\left(  \left\vert \log D_{2}h_{j,n_{k}%
}\right\vert \vee\log\log n_{k}\right)  },\hspace{5mm}j\geq0.
$$
Applying Inequality A.1 in the Appendix with
$$
M=\kappa\quad \textrm{ and }\quad \sigma_{\mathcal{G}}^{2}=\sigma_{\mathcal{F}_{j,n_{k}}%
}^{2}\leq D_{0}h_{j,n_{k}},
$$
we get for any $t>0,$%
$$
I\!\!P\left\{  \max_{n_{k-1}\leq n\leq n_{k}}||\sqrt{n}\alpha_{n}%
||_{\mathcal{F}_{j,n_{k}}}\geq A_{1}(D_{1}c_{j,k}+t)\right\}
$$%
$$
\leq2\left[  \exp\left(  -A_{2}t^{2}/\left(  D_{0}n_{k}h_{j,n_{k}}\right)
\right)  +\exp(-A_{2}t/\kappa)\right]  .
$$
Set for any $\rho>1,$ $j\geq0$ and $k\geq1,$%
$$
p_{j,k}(\rho):=I\!\!P\left\{  \max_{n_{k-1}\leq n\leq n_{k}}||\sqrt{n}%
\alpha_{n}||_{\mathcal{F}_{j,n_{k}}}\geq A_{1}\left(  D_{1}+\rho\right)
c_{j,k}\right\}  .
$$
As we have $c_{j,k}/\sqrt{n_{k}h_{j,n_{k}}}\geq\sqrt{\log\log n_{k}}$, we
readily obtain for $j\geq0,$
$$
p_{j,k}(\rho)\leq2\left[  \exp\left(  -\mbox{ }\frac{\rho^{2}A_{2}}{D_{0}}%
\log\log n_{k}\right)  +\exp\left(  -\frac{\sqrt{c}\rho A_{2}}{\kappa}%
\sqrt{\log n_{k}\log\log n_{k}}\right)  \right]  ,
$$
which for $\gamma=$ $\frac{A_{2}}{D_{0}}\wedge\frac{\sqrt{c}A_{2}}{\kappa}$
implies
$$
p_{j,k}(\rho)\leq4\exp\left(  -\mbox{ }\rho\gamma\log\log n_{k}\right)  .
$$
Thus
$$
P_{k}(\rho):=\sum_{j=0}^{l_{n_{k}}-1}p_{j,k}(\rho)\leq4l_{n_{k}}\left(  \log
n_{k}\right)  ^{-\rho\gamma},
$$
which by (\ref{LL}), for all large $k$ and large enough $\rho>1$
$$
P_{k}(\rho)\leq8\left(  \log n_{k}\right)  ^{1-\rho\gamma}=8\left(  \frac
{1}{k\log2}\right)  ^{\rho\gamma-1}\leq k^{-2}.
$$
Notice that by definition of $l_{n}$, for large $k$
$$
2h_{l_{n_{k}},n_{k}}=h_{l_{n_{k}}+1,n_{k}}\geq2h_{0},
$$
which implies that we have for $n_{k-1}\leq n\leq n_{k}$
$$
\left[  \frac{c\log n}{n},h_{0}\right]  \subset\left[  \frac{c\log n_{k}%
}{n_{k}},h_{l_{n_{k}},n_{k}}\right]  .
$$
Thus for all large enough $k$ and $n_{k-1}\leq n\leq n_{k},$%
$$
A_{k}(\rho):=
$$%
$$
\left\{  \max_{n_{k-1}\leq n\leq n_{k}}\sup_{g\in\mathcal{G}}\sup_{\frac{c\log
n}{n}\leq h\leq h_{0}}\frac{\sqrt{nh}\Vert g_{n,h}-\mathrm{I\!E}g_{n,h}\Vert_{\infty
}}{\sqrt{|\log h|\vee\log\log n}}>2A_{1}(D_{1}+\rho)\right\}
$$%
$$
\subset\bigcup_{j=0}^{l_{n_{k}}-1}\left\{  \max_{n_{k-1}\leq n\leq n_{k}%
}||\sqrt{n}\alpha_{n}||_{\mathcal{F}_{j,n_{k}}}\geq A_{1}(D_{1}+\rho
)c_{j,k}\right\}  .
$$
It follows now for large enough $\rho$ that
$$
I\!\!P\left\{  A_{k}(\rho)\right\}  \leq P_{k}(\rho)\leq k^{-2},
$$
which by the Borel-Cantelli lemma implies our theorem. $\Box$

\section{Application to local polynomial regression function estimators}

\textit{In this section we shall always assume that the assumptions of Theorem
1 hold (in particular, that }$K$\textit{ has support contained in} $\left[
-1/2,1/2\right]  $) \textit{and }$I$ \textit{is a fixed compact interval in
}$\mathrm{I\!R}$. \textit{We shall also assume that} $K\geq0$. \smallskip

\subsection{Estimating the regression function by local polynomials}

Let $(X,Y),(X_{1},Y_{1}),\ldots,(X_{n},Y_{n})$ be i.i.d. $2$-dimensional
random vectors and write
$$
g(x):=\mathrm{I\!E}\left[  Y||X=x\right]
$$
for the regression function. Suppose that $g(x)$ is $(p+1)$ times
differentiable on $J=I^{\eta}$, then we can approximate $g(x)$ locally around
$x_{0}\in I$ by a polynomial of order $p$ (Taylor):
$$
g(x)\approx g(x_{0})+g^{\prime}(x_{0})(x-x_{0})+\ldots+\frac{g^{(p)}(x_{0}%
)}{p!}(x-x_{0})^{p}.
$$
Then consider the weighted least-squares regression problem (WLS)
\begin{equation}
{\text{argmin}}_{\pmb{\beta}\in I\!\!R^{p+1}}\;\;\frac{1}{nh}\sum_{i=1}%
^{n}[Y_{i}-\sum_{j=0}^{p}\beta_{j}(X_{i}-x_{0})^{j}]^{2}K\left(  \frac
{x_{0}-X_{i}}{h}\right)  . \label{WLS}%
\end{equation}
It is clear that if $\pmb{\hat\beta}\in I\!\!R^{p+1}$ is the solution of the
WLS problem in (\ref{WLS}), we obtain an estimator $\hat{g}_{n,h}^{(p)}%
(x_{0})$ of $g(x_{0})$ by taking it be $\hat{\beta}_{0}$, the first component
of $\pmb{\hat\beta}$. At the same time we obtain estimators of the derivatives
of the regression function up to order $p$. To solve (\ref{WLS}), first note
that it can be written in a matrix notation:
\begin{equation}
{\text{argmin}}_{\pmb\beta\in I\!\!R^{p+1}}\left(  \pmb Y-\pmb X_{x_{0}%
}\pmb\beta\right)  ^{t}\pmb W_{x_{0}}\left(  \pmb Y-\pmb X_{x_{0}}%
\pmb\beta\right)  , \label{WLS matrix}%
\end{equation}
where $\pmb W_{x_{0}}=\left(  nh\right)  ^{-1}diag\left(K\left(
\frac{x_{0}-X_{i}}{h}\right)  \right)  \in I\!\!R^{n\times n}$, and $\pmb
X_{x_{0}}\in I\!\!R^{n\times(p+1)},\pmb Y\in I\!\!R^{n\times1}$ and
$\pmb\beta\in I\!\!R^{(p+1)\times1}$ are defined as%

$$
\pmb X_{x_{0}}:=\left(
\begin{array}
[c]{cccc}%
1 & (X_{1}-x_{0}) & \cdots & (X_{1}-x_{0})^{p}\\
\vdots & \vdots &  & \vdots\\
1 & (X_{n}-x_{0}) & \cdots & (X_{n}-x_{0})^{p}%
\end{array}
\right)  ,\text{ }\pmb Y:=\left(
\begin{array}
[c]{c}%
Y_{1}\\
\vdots\\
Y_{n}%
\end{array}
\right)  ,\text{ }\pmb\beta:=\left(
\begin{array}
[c]{c}%
\beta_{0}\\
\vdots\\
\beta_{p}%
\end{array}
\right)  .
$$
If we set
$$
L(x_{0}):=\frac1{nh}\sum_{i=1}^{n}[Y_{i}-\sum_{j=0}^{p}\beta_{j}(X_{i}%
-x_{0})^{j}]^{2}K\left(  \frac{ x_{0}-X_{i}}h\right)  ,
$$
it is not too difficult to see that for $k=0,\ldots,p$, the partial
derivatives can be written as
$$
\frac{\partial L(x_{0})}{\partial\beta_{k}}=-2(\pmb Y-\pmb X_{x_{0}}%
\pmb\beta)^{t}\pmb W_{x_{0}}\pmb X_{x_{0}}\pmb e_{k}^{t},
$$
where $\pmb e_{k}$ is the $k$-th unit vector in $I\!\!R^{p+1}$. So by setting
the partial derivatives equal to zero, we obtain that the solution
$\pmb{\hat \beta}$ of the WLS problem (\ref{WLS}) must satisfy
$$
\pmb Y^{t}\pmb W_{x_{0}}\pmb X_{x_{0}}=\pmb{\hat \beta}^{t}\pmb X_{x_{0}}%
^{t}\pmb W_{x_{0}}\pmb X_{x_{0}}.
$$
Assuming that
$$
\mathcal{S}_{x_{0}}:=\pmb X_{x_{0}}^{t}\pmb W_{x_{0}}\pmb X_{x_{0}},
$$
is invertible, we can compute the solution by
$$
\pmb{\hat \beta}_{x_{0}}=\left(  \pmb X_{x_{0}}^{t}\pmb W_{x_{0}}\pmb
X_{x_{0}}\right)  ^{-1}\pmb X_{x_{0}}^{t}\pmb W_{x_{0}}\pmb Y.
$$
We shall show that \emph{asymptotically }the inverse matrix of $\mathcal{S}%
_{x_{0}}$ always exists. To see this, consider for $0\leq j\leq2p$ the
functions
$$
H^{(j)}(u):=\left(  -u\right)  ^{j}K(u).
$$
Since we assume $K$ to be bounded with support contained in $\left[
-1/2,1/2\right]  $, we see that each $H^{(j)}\in L_{1}(I\!\!R)$ and has
support contained in $\left[  -1/2,1/2\right]  $. Now for each $j\geq0$ define
the bounded function%

$$
\phi_{j}\left(  u\right)  =\left(  -u\right)  ^{j}1\left\{  u\in\left[
-1/2,1/2\right]  \right\}  .
$$
Since this function is of bounded variation, the class
$$
\left\{  \phi_{j}((x-\cdot)/h):h>0,x\in\mathrm{I\!R}\right\}
$$
satisfies (A.2). (See Lemma 22 of \cite{8a}.) Thus the class
$\mathcal{K}$, as defined in (\ref{KK}) is assumed to be pointwise measurable
and satisfies (A.2). By Lemma A.1 in the\ Appendix, for each $j=0,\dots,2p$, the
class
$$
\mathcal{G}_{j}:=\left\{  H^{(j)}((x-\cdot)/h):h>0,x\in\mathrm{I\!R}\right\}
$$
also fulfills (A.2). Moreover, it is easily checked that each $\mathcal{G}%
_{j}$ is pointwise measurable. Hence the assumptions of Theorem 2 hold and we
can infer that for each $0\leq j\leq2p$, and sequence $a_{n}$ satisfying%
\begin{equation}
a_{n}\searrow0\text{ and }na_{n}/\log n\rightarrow\infty, \label{aa}%
\end{equation}
we have
\begin{equation}
\sup_{x_{0}\in I}\sup_{a_{n}\leq h\leq h_{0}}\left\vert H_{n,h}^{(j)}%
(x_{0})-\mathrm{I\!E}H_{n,h}^{(j)}(x_{0})\right\vert \longrightarrow0,\hspace
{5mm}a.s., \label{aa1}%
\end{equation}
where
$$
H_{n,h}^{(j)}(x_{0}):=\frac{1}{nh}\sum_{i=1}^{n}H^{(j)}\left(  \frac
{x_{0}-X_{i}}{h}\right)  .
$$
Notice that
$$
\mathrm{I\!E}H_{n,h}^{(j)}(x_{0})=\frac{1}{h}\int_{\mathrm{I\!R}}H^{(j)}\left(
\frac{x_{0}-t}{h}\right)  f(t)dt=:f* H_{h}^{(j)}(x_{0}),
$$
and since $f$ is continuous on $J=I^{\eta}$ with $I$ being a compact interval,
we can use Lemma A.2 in the Appendix to get that as $h\searrow0$,
\begin{equation}
\sup_{x_{0}\in I}\left\vert \mathrm{I\!E}H_{n,h}^{(j)}(x_{0})-f(x_{0})\int
_{\mathrm{I\!R}}\left(  -u\right)  ^{j}K(u)du\right\vert \longrightarrow0.
\label{aa2}%
\end{equation}
Hence, it follows immediately by $(\ref{aa1})$ and $(\ref{aa2})$, that
uniformly in $x_{0}\in I$ and for $a_{n}<b_{n}$ with $a_{n}$ satisfying
$(\ref{aa})$ and $b_{n}\searrow0,$
\begin{equation}
\sup_{a_{n}\leq h\leq b_{n}}\left\vert H_{n,h}^{(j)}(x_{0})-f(x_{0}%
)\int_{\mathrm{I\!R}}\left(  -u\right)  ^{j}K(u)du\right\vert \longrightarrow
0,\hspace{5mm}a.s. \label{conv}%
\end{equation}
Next consider the Hilbert space $\mathcal{L}(\mathrm{I\!R},Kd\lambda)$
consisting of all the measurable functions $\phi:\mathrm{I\!R}\rightarrow
\mathrm{I\!R}$ such that
$$
\int_{\mathrm{I\!R}}\phi^{2}(u)K(u)du<\infty.
$$
As usual, $\phi_{1}=\phi_{2}$ if $\int_{\mathrm{I\!R}}(\phi_{1}-\phi_{2}%
)^{2}(u)K(u)du=0$; that is, each $\phi\in\mathcal{L}(\mathrm{I\!R},Kd\lambda)$
represents an equivalence class of functions. Now let
$$
G:=\left(  \int_{\mathrm{I\!R}}\left(  -u\right)  ^{j+k}K(u)du\right)
_{j=0,k=0}^{p,p},
$$
then $G$ is the Gramian matrix of the set of functions $\{\varphi_{j}%
:\varphi_{j}(x)=\left(  -x\right)  ^{j},$ $j=0,\dots,p\}$ and these functions
belong to $\mathcal{L}(\mathrm{I\!R},Kd\lambda)$ since $K$ has compact
support. It is known that $G$ is nonsingular if the functions are linearly
independent. Hence, in our case, $G$ will always be invertible. (Here we use
$K\geq0$ and $0<\int_{\mathrm{I\!R}}K(u)du<\infty$.) To see that
$\mathcal{S}_{x_{0}}$ is invertible as well, recall that the function
$M\rightarrow\det M$ with $M\in\mathcal{M}_{p+1}(\mathrm{I\!R})$ is
continuous, and that by $(\ref{aa1})$ and $(\ref{aa2})$, with probability one,
the components of
$$
A_{x_{0}}:=\left(  H_{n,h}^{(j+k)}(x_{0})\right)  _{j=0,k=0}^{p,p},
$$
converge uniformly in $x_{0}\in I$ and $a_{n}\leq h\leq b_{n}$ with
$b_{n}\searrow0$ to those of $f(x_{0})G$. Hence, since we assume $f$ to be
strictly positive on $J=I^{\eta}$, for $n$ large enough, uniformly in
$x_{0}\in I$, we have $\det A_{x_{0}}>0$. Now let $\mathcal{H}_{p}%
:=diag\{1,h,\dots,h^{p}\}$, note that
$$
\mathcal{S}_{x_{0}}=\mathcal{H}_{p}A_{x_{0}}\mathcal{H}_{p},
$$
and observe that $\det\mathcal{S}_{x_{0}}=h^{p(p+1)}\det A_{x_{0}}$, so for
$n$ large enough, uniformly in $x_{0}\in I$ and $a_{n}\leq h\leq b_{n}$,
$\mathcal{S}_{x_{0}}$ will have a positive determinant, showing that
\emph{asymptotically}, $\mathcal{S}_{x_{0}}$ is nonsingular and
invertible.\newline

From the above it follows that with probability one, for all large $n$,
uniformly in $x_{0}\in I$ and $a_{n}\leq h\leq b_{n}$, the local polynomial
regression estimator of $g(x_{0})$ is given by
$$
\hat{g}_{n,h}^{(p)}(x_{0})=\pmb e_{1}\mathcal{S}_{x_{0}}^{-1}\pmb X_{x_{0}%
}^{t}\pmb W_{x_{0}}\pmb Y.
$$
The difficulty is to determine $\mathcal{S}_{x_{0}}^{-1}$ explicitly,
especially when $p$ becomes large. Moreover, it is not possible to find a nice
general formula for $\hat{g}_{n,h}^{(p)}(x_{0})$, since the calculation of
$\mathcal{S}_{x_{0}}^{-1}$ and $\hat{g}_{n,h}^{(p)}(x_{0})$ becomes more
complex as $p$ increases. However, we shall see in the next section that
$\hat{g}_{n,h}^{(p)}(x_{0})$ can be easily computed for $p=0,1,2$.

\subsection{Uniform in bandwidth consistency}

We shall now discuss uniform in bandwidth consistency of $\hat{g}_{n,h}^{(p)}$
on a compact interval $I$. Define the functions
\begin{align*}
&  \tilde{f}_{n,h,j}(x):=\frac{1}{nh}\sum_{i=1}^{n}\left(  \frac{X_{i}-x}%
{h}\right)  ^{j}K\left(  \frac{x-X_{i}}{h}\right)  ,\;\;j=0,\ldots,2p,\\
&  \tilde{r}_{n,h,j}(x):=\frac{1}{nh}\sum_{i=1}^{n}Y_{i}\left(  \frac{X_{i}%
-x}{h}\right)  ^{j}K\left(  \frac{x-X_{i}}{h}\right)  ,\;\;j=0,\ldots,p.
\end{align*}
By Theorem 2,
$$
\limsup_{n\rightarrow\infty}\sup_{\frac{c\log n}{n}\leq h\leq h_{0}}%
\max_{0\leq j\leq 2p}\frac{\sqrt{nh}\left\Vert \tilde{f}_{n,h,j}-\mathrm{I\!E}%
\tilde{f}_{n,h,j}\right\Vert _{I}}{\sqrt{ |\log h|\vee\log\log n }}%
<\infty,\text{ }a.s.
$$
and by Theorem 1 with obvious identifications and $K$ replaced by $H^{(j)}$,
$$
\limsup_{n\rightarrow\infty}\sup_{c(\log n/n)^{\gamma}\leq h\leq h_{0}}%
\max_{0\leq j\leq p}\frac{\left\Vert \tilde{r}_{n,h,j}-\mathrm{I\!E}\tilde
{r}_{n,h,j}\right\Vert _{I}}{\sqrt{nh\left(  |\log h| \vee\log\log n\right)
}}<\infty,\text{ }a.s.
$$
For $j\geq0$, set
$$
\mu_{j}:=\int_{\mathrm{I\!R}}(-u)^{j}K(u)du,
$$
and define
\begin{align*}
&  f_{j}(x):=\mu_{j}f_{X}(x),\;\;\;j=0,\ldots,2p,\\
&  r_{j}(x):=\mu_{j}\int_{\mathrm{I\!R}}yf(x,y)dy,\;\;\;j=0,\ldots,p.
\end{align*}
Lemma A.2 gives (also see $(\ref{aa2})$) that for all $0\leq j\leq2p$,
\begin{equation}
\sup_{a_{n}\leq h\leq b_{n}}\Vert\mathrm{I\!E}\tilde{f}_{n,h,j}-f_{j}%
\Vert_{\infty}\longrightarrow0. \label{ff}%
\end{equation}
Now define the function
$$
\varphi(x):=\int_{\mathrm{I\!R}}yf\left(  x,y\right)  dy,\;\;\; x\in J,
$$
and introduce the assumption:
\begin{equation}
\mbox{for all }\mbox{ }x\in J,\mbox{ }\lim_{x^{\prime}\rightarrow
x}f(x^{\prime},y)=f\left(  x,y\right)  \mbox{ for almost every }y\in
\mathrm{I\!R}\mathbf{.} \label{f2}%
\end{equation}
Then by an argument based on the Lebesgue dominated convergence theorem, using
assumptions $(\ref{f3})$ along with $(\ref{mm})$ or $(\ref{PP})$, one readily
shows that $\varphi$ is bounded and continuous on $J$. Applying Lemma A.2, we
get that for all $0\leq j\leq p$,
\begin{equation}
\sup_{a_{n}\leq h\leq b_{n}}\Vert\mathrm{I\!E}\tilde{r}_{n,h,j}-r_{j}\Vert
_{I}\longrightarrow0. \label{rr}%
\end{equation}

From these observations, we easily conclude that for all smooth functions
$\Phi:I\!\!R^{3p+2}\rightarrow I\!\!R$ and suitable sequences $0<a_{n}<b_{n}$
depending on Theorem 1 and whether $(\ref{mm})$ or $(\ref{PP})$ holds, with
probability $1$,
\begin{equation}
\sup_{a_{n}\leq h\leq b_{n}}\left\Vert \Phi\left(  \tilde{f}_{n,h,0}%
,\ldots,\tilde{f}_{n,h,2p},\tilde{r}_{n,h,0},\ldots,\tilde{r}_{n,h,p}\right)
-\Phi\left(  f_{0},\ldots,f_{2p},r_{0},\ldots,r_{p}\right)  \right\Vert
_{I}\longrightarrow0. \label{obs2}%
\end{equation}
When $(\ref{mm})$ is in force, we assume that $a_{n}$ satisfies $(\ref{aa})$,
and when $(\ref{PP})$ holds that $a_{n}=c(\log n/n)^{\gamma}$ for $\gamma>1$.

\paragraph{Calculation for $p=0$.}

In this case we get the usual Nadaraya-Watson regression estimator:
$$
\hat{g}_{n,h}^{(0)}(x_{0})=\frac{\sum_{i=1}^{n}Y_{i}K\left(  \frac{x_{0}%
-X_{i}}{h}\right)  }{\sum_{i=1}^{n}K\left(  \frac{x_{0}-X_{i}}{h}\right)  }=\frac{\tilde{r}_{n,h,0}(x_{0})}{\tilde{f}_{n,h,0}(x_{0})}.
$$
So applying (\ref{obs2}) with $\Phi(x_{1},x_{2})=x_{2}/x_{1}$, we get that
uniformly in $x_{0}\in I$,
$$
\sup_{a_{n}\leq h\leq b_{n}}\left\Vert \hat{g}_{n,h}^{(0)}-g\right\Vert
_{I}\longrightarrow0,\;\;a.s.,
$$
proving the uniform in bandwidth consistency of the Nadaraya-Watson
estimator.\newline

\textit{From now on, for ease of notation we shall omit the subscripts }%
$x_{0}$\textit{, as well as the argument }$(x_{0})$\textit{ in all the
functions that we defined above.}

\paragraph{Calculation for $p=1$.}

This is the local linear regression estimator, where $\mathcal{S}$ and $\pmb
X^{t}\pmb W\pmb Y$ are given by
$$
\mathcal{S}=\left(
\begin{array}
[c]{cc}%
nh\tilde{f}_{n,h,0} & nh^{2}\tilde{f}_{n,h,1}\\
nh^{2}\tilde{f}_{n,h,1} & nh^{3}\tilde{f}_{n,h,2}%
\end{array}
\right)  ,\hspace{1cm}\pmb X^{t}\pmb W\pmb Y=\left(
\begin{array}
[c]{c}%
nh\tilde{r}_{n,h,0}\\
nh^{2}\tilde{r}_{n,h,1}%
\end{array}
\right)  ,
$$
such that
$$
\mathcal{S}^{-1}\pmb X^{t}\pmb W\pmb Y=\frac{1}{\tilde{f}_{n,h,0}\tilde
{f}_{n,h,2}-\tilde{f}_{n,h,1}^{2}}\left(
\begin{array}
[c]{c}%
\tilde{f}_{n,h,2}\tilde{r}_{n,h,0}-\tilde{f}_{n,h,1}\tilde{r}_{n,h,1}\\
\tilde{f}_{n,h,0}\tilde{r}_{n,h,1}-\tilde{f}_{n,h,1}\tilde{r}_{n,h,0}%
\end{array}
\right)  .
$$
Hence, the local linear estimator of the regression function is given by
$$
\hat{g}_{n,h}^{(1)}=\frac{\tilde{f}_{n,h,2}\tilde{r}_{n,h,0}-\tilde{f}%
_{n,h,1}\tilde{r}_{n,h,1}}{\tilde{f}_{n,h,0}\tilde{f}_{n,h,2}-\tilde
{f}_{n,h,1}^{2}}.
$$
So applying (\ref{obs2}) with $\Phi(x_{1},\ldots,x_{5})=\frac{x_{3}x_{4}%
-x_{2}x_{5}}{x_{1}x_{3}-x_{2}^{2}}$, we obtain after a little algebra based on
the definitions of $f_{j}$ and $r_{j}$, the uniform in bandwidth consistency of
this local linear estimator:
$$
\sup_{a_{n}\leq h\leq b_{n}}\left\Vert \hat{g}_{n,h}^{(1)}-g\right\Vert
_{I}\longrightarrow0,\;\;a.s.
$$

\paragraph{Calculation for $p=2$.}

As we have seen in the case $p=1$, the main work in deriving $\hat{g}%
_{n,h}^{(2)}$ is to determine $\mathcal{S}^{-1}$. Now $\mathcal{S}$ is a
$3\times3$-matrix, so we can still write down the inverse without
difficulties. After some calculations, we obtain (disregarding $nh^{j}$ factors):%

$$
\mathcal{S}^{-1}=
$$%
$$
\frac{1}{\det\mathcal{S}}\left(
\begin{array}
[c]{lcr}%
\tilde{f}_{n,h,2}\tilde{f}_{n,h,4}-\tilde{f}_{n,h,3}^{2} & \tilde{f}%
_{n,h,2}\tilde{f}_{n,h,3}-\tilde{f}_{n,h,1}\tilde{f}_{n,h,4} & \tilde
{f}_{n,h,1}\tilde{f}_{n,h,3}-\tilde{f}_{n,h,2}^{2}\\
\tilde{f}_{n,h,2}\tilde{f}_{n,h,3}-\tilde{f}_{n,h,1}\tilde{f}_{n,h,4} &
\tilde{f}_{n,h,0}\tilde{f}_{n,h,4}-\tilde{f}_{n,h,2}^{2} & \tilde{f}%
_{n,h,1}\tilde{f}_{n,h,2}-\tilde{f}_{n,h,0}\tilde{f}_{n,h,3}\\
\tilde{f}_{n,h,1}\tilde{f}_{n,h,3}-\tilde{f}_{n,h,2}^{2} & \tilde{f}%
_{n,h,1}\tilde{f}_{n,h,2}-\tilde{f}_{n,h,0}\tilde{f}_{n,h,3} & \tilde
{f}_{n,h,0}\tilde{f}_{n,h,2}-\tilde{f}_{n,h,1}^{2}\\
&  &
\end{array}
\right)  ,
$$
and%
$$
\pmb X^{t}\pmb W\pmb Y=\left(
\begin{array}
[c]{c}%
\tilde{r}_{n,h,0}\\
\widetilde{r}_{n,h,1}\\
\widetilde{r}_{n,h,2}%
\end{array}
\right)  ,
$$
eventually yielding%

$$
\hat{g}_{n,h}^{(2)}=
$$%
$$
\frac{(\tilde{f}_{n,h,2}\tilde{f}_{n,h,4}-\tilde{f}_{n,h,3}^{2})\tilde
{r}_{n,h,0}+(\tilde{f}_{n,h,2}\tilde{f}_{n,h,3}-\tilde{f}_{n,h,1}\tilde
{f}_{n,h,4})\tilde{r}_{n,h,1}+(\tilde{f}_{n,h,1}\tilde{f}_{n,h,3}-\tilde
{f}_{n,h,2}^{2})\tilde{r}_{n,h,2}}{\tilde{f}_{n,h,0}\tilde{f}_{n,h,2}\tilde
{f}_{n,h,4}-\tilde{f}_{n,h,0}\tilde{f}_{n,h,3}^{2}-\tilde{f}_{n,h,1}^{2}%
\tilde{f}_{n,h,4}+2\tilde{f}_{n,h,1}\tilde{f}_{n,h,2}\tilde{f}_{n,h,3}%
-\tilde{f}_{n,h,2}^{3}}.
$$
So using the function
$$
\Phi(x_{1},\ldots,x_{8})=\frac{(x_{3}x_{5}-x_{4}^{2})x_{6}+(x_{3}x_{4}%
-x_{2}x_{5})x_{7}+(x_{2}x_{4}-x_{3}^{2})x_{8}}{x_{1}x_{3}x_{5}-x_{1}x_{4}%
^{2}-x_{2}^{2}x_{5}+2x_{2}x_{3}x_{4}-x_{3}^{3}}
$$
in (\ref{obs2}), we infer after some algebra based on the definitions of
$f_{j}$ and $r_{j}$, the uniform in bandwidth consistency of this local
quadratic regression function estimator.

\paragraph{Calculation for larger $p$.}

In principle it is possible to write down an explicit formula for the local
polynomial estimator $\hat{g}_{n}^{(p)}(x_{0})$ for any $p\geq0$, by first
computing the inverse of $\mathcal{S}_{x_{0}}$, multiplying it by $\pmb
X_{x_{0}}^{t}\pmb W_{x_{0}}\pmb Y$ and then by taking the first component of
the resulting vector. But the difficulty lies in determining $\mathcal{S}%
_{x_{0}}^{-1}$.

\paragraph{Remark.}

It was pointed in \cite{4a} and \cite{1a} that
these methods can be used to study the uniform in bandwidth consistency of
local polynomial regression estimators.

\section{Appendix}

Let $X,$ $X_{1},\ldots,X_{n}$ be i.i.d. from a probability space
$(\mathcal{X},\mathcal{A},P)$ with common distribution $\mu$. Let
$\mathcal{G}$ be a \textit{pointwise measurable class} of real valued
functions defined on $\mathcal{X}$, i.e. we assume that there exists a
countable subclass $\mathcal{G}_{0}$ of $\mathcal{G}$ so that we can find for
any function $g$ in $\mathcal{G}$ a sequence of functions $\{g_{m}\}$ in
$\mathcal{G}_{0}$ for which $g_{m}(x)\rightarrow g(x),$ $x\in\mathcal{X}$.
(See Example 2.3.4, \cite{17}.) Further let
$\varepsilon_{1},\ldots,\varepsilon_{n}$ be a sequence of independent
Rademacher random variables independent of $X_{1},\ldots,X_{n}$.

The following inequality is essentially due to \cite{14} (see \cite{4}). \medskip

\noindent\textbf{Inequality A.1} \textit{Let} $\mathcal{G}$ \textit{be a
pointwise measurable class of functions satisfying for some} $0<M<\infty$
$$
||g||_{\infty}\leq M,\mbox{ }g\in\mathcal{G},
$$
\textit{then} \textit{for all} $t>0$ \textit{we have for suitable finite
constants} $A_{1},A_{2}>0$,
$$
\mathrm{I\!P}\left\{  \max_{1\leq m\leq n}||\sqrt{m}\alpha_{m}||_{\mathcal{G}}\geq
A_{1}(\mathrm{I\!E}||\sum_{i=1}^{n}\varepsilon_{i}g(X_{i})||_{\mathcal{G}}+t)\right\}
$$
$$
\leq2(\exp(-A_{2}t^{2}/n\sigma_{\mathcal{G}}^{2})+\exp(-A_{2}t/M)),
$$
\textit{where} $\sigma_{\mathcal{G}}^{2}=\sup_{g\in\mathcal{G}}Var(g(X))$%
.\medskip

It enables us to reduce many problems on almost sure convergence to
investigating the moment quantity
$$
\mu_{n}:=\mathrm{I\!E}||\sum_{i=1}^{n}\varepsilon_{i}g(X_{i})||_{\mathcal{G}}.
$$

The following proposition proved in \cite{4} is very helpful
for obtaining bounds on this quantity, when the class $\mathcal{G}$ has a
polynomial covering number. Let $G$ be a finite valued measurable function
satisfying for all $x\in\mathcal{X}$
$$
G(x)\geq\sup_{g\in\mathcal{G}}|g(x)|,
$$
and define
$$
N(\epsilon,\mathcal{G}):=\sup_{Q}N(\epsilon\sqrt{Q(G^{2})},\mathcal{G},d_{Q}),
$$
where the supremum is taken over all probability measures $Q$ on
$(\mathcal{X},\mathcal{A})$ for which $0<Q(G^{2})<\infty$ and $d_{Q}$ is the
$L_{2}(Q)$--metric. As usual $N(\epsilon,\mathcal{G},d)$ is the minimal number
of balls $\{g:d(g,f)<\epsilon\}$ of $d$-radius $\epsilon$ needed to cover
$\mathcal{G}$.\medskip

\noindent\textbf{Proposition A.1} \textit{Let} $\mathcal{G}$ \textit{be a
pointwise measurable class of bounded functions such that for some constants}
$\beta,\nu,C>1,$ $\sigma\leq1/(8C)$ \textit{and function} $G$ \textit{as
above}, \textit{the following four conditions hold:}
$$
\mathrm{I\!E}[G^{2}(X)]\leq\beta^{2};\eqno(A.1)
$$%
$$
N(\epsilon,\mathcal{G})\leq C\epsilon^{-\nu},\mbox{ }0<\epsilon<1;\eqno(A.2)
$$%
$$
\sigma_{0}^{2}:=\sup_{g\in\mathcal{G}}\mathrm{I\!E}[g^{2}(X)]\leq\sigma^{2};\eqno(A.3)
$$%
$$
\sup_{g\in\mathcal{G}}||g||_{\infty}\leq\frac{1}{2\sqrt{\nu+1}}\sqrt
{n\sigma^{2}/\log(\beta\vee1/\sigma)}.\eqno(A.4)
$$
\textit{Then we have} \textit{for a universal constant }$A$%
$$
\mathrm{I\!E}||\sum_{i=1}^{n}\varepsilon_{i}g(X_{i})||_{\mathcal{G}}\leq A\sqrt{\nu
n\sigma^{2}\log(\beta\vee1/\sigma)}.\eqno(A.5)
$$
Another version of Proposition A.1 has been proved by \cite{6}. For refinements, consult \cite{4a} and \cite{8}.\medskip

We shall also require the following two lemmas. The first is proved in \cite{4}.\medskip

Here is Lemma A.1 of \cite{4}.\medskip

\noindent\textbf{Lemma A.1} \textit{Let} $\mathcal{F}$ \textit{and}
$\mathcal{G}$ \textit{be two classes of real valued measurable functions on}
$\mathcal{X}$ \textit{satisfying}
$$
|f(x)|\leq F(x),\mbox{ }f\in\mathcal{F},x\in\mathcal{X}%
$$
\textit{where } $F$ \textit{is a finite valued measurable envelope function
on} $\mathcal{X}$;
$$
\parallel g\parallel_{\infty}\leq M,\mbox{ }g\in\mathcal{G},
$$
\textit{where} $M>0$ \textit{is a finite constant. Assume that for all}
$p$-\textit{measures }$Q$ \textit{with} $0<Q(F^{2})<\infty,$
$$
N(\epsilon\sqrt{Q(F^{2})},\mathcal{F},d_{Q})\leq C_{1}\epsilon^{-\nu_{1}%
},0<\epsilon<1,
$$
\textit{and for all }$p$-\textit{measures} $Q$,
$$
N(\epsilon M,\mathcal{G},d_{Q})\leq C_{2}\epsilon^{-\nu_{2}},0<\epsilon<1,
$$
\textit{where} $\nu_{1},\nu_{2},C_{1},C_{2}\geq1$ \textit{are}
\textit{suitable constants. Then we have for all} $p$-\textit{measures} $Q$,
\textit{with }$Q(F^{2})<\infty$,
$$
N(\epsilon M\sqrt{Q(F^{2})},\mathcal{F}\mathcal{G},d_{Q})\leq C_{3}%
\epsilon^{-\nu_{1}-\nu_{2}},0<\epsilon<1,
$$
\textit{for some finite constant} $0<C_{3}<\infty.$\medskip

The next lemma can be inferred from results in \cite[pp. ~62--65]{9b}.\textit{\smallskip}

\noindent\textbf{Lemma A.2. }\textit{Let }$\varphi$\textit{\ be a measurable
function on }$\mathrm{I\!R}^{d}$\textit{, which for some }$\gamma>0$
\textit{is bounded and uniformly continuous on }$D_{\gamma},$ \textit{where
}$D$ \textit{is a closed subset of }$\mathrm{I\!R}^{d}\mathit{\ }%
$\textit{and}
$$
D_{\gamma}=\left\{  x\in\mathrm{I\!R}^{d}:|x-y|\leq\gamma,y\in D\right\}  .
$$
\textit{\ Then for any }$L_{1}(\mathrm{I\!R}^{d})$ \textit{function }%
$H$\textit{, which is equal to zero for }$x\notin I^{d}$%
$$
\sup_{z\in D}|\varphi\ast H_{h}(z)-I(H)\varphi(z)|\rightarrow
0,\mbox{{\it \ as} }h\searrow0,
$$
\textit{where} $I(H)=\int_{\mathrm{I\!R}^{d}}H(u)du$ \textit{and }$\varphi\ast
H_{h}(z):=h^{-1}\int_{\mathrm{I\!R}^{d}}\varphi(x)H\left(  h^{-1/d}\left(
z-x\right)  \right)  dx.$

\paragraph{Acknowledgements.}

David Mason's research was partially supported by an NSF Grant. Julia Dony's
research is financed by a PhD grant from the Institute for the Promotion of
Innovation through Science and Technology in Flanders (IWT Vlaanderen).

\bibliographystyle{plain}
\bibliography{mikfinall}


\end{document}